\documentclass[12pt]{amsart}
\usepackage{amsmath,amssymb,amsthm}
\usepackage{graphicx}

\theoremstyle{definition}

\newcommand{\ds}{\displaystyle}
\begin{document}

\title{ISOGONAL CONJUGACY AND FERMAT PROBLEMS}
\author{Georgi Ganchev and Nikolai Nikolov}
\address{Bulgarian Academy of Sciences, Institute of Mathematics and Informatics,
Acad. G. Bonchev Str. bl. 8, 1113 Sofia, Bulgaria}
\email{ganchev@math.bas.bg}
\email{nik@math.bas.bg}

\begin{abstract}
We consider three types of isogonal conjugacy of two points with
respect to a given triangle and characterize any of these types by a geometric equality.
As an application to the Fermat problem with positive weights, we prove that in the
general case the given weights determine uniquely a point X and the solution to the
Fermat problem is the point Y, which is isogonally conjugate of type I to the point X.
We obtain a similar characterization of the solution to the Fermat problem in the case
of mixed weights as well.
\end{abstract}
\maketitle

\section{Introduction}

In the present paper we analyze the following problems:
\vskip 2mm
{\it Given a $\triangle ABC$ and two points $X$, $Y$ in its plane. Prove that the
equality
$$\frac{AX\,AY}{AB\,AC}+\frac{BX\,BY}{BA\,BC}+\frac{CX\,CY}{CA\,CB}=1$$
holds   if and only if $X$ and $Y$ are isogonally conjugate of type
$I$}.
\vskip 4mm
{\it Given a $\triangle ABC$ and two points $X$, $Y$ in its plane. Prove that the
equality
$$-\frac{AX\,AY}{AB\,AC}+\frac{BX\,BY}{BA\,BC}+\frac{CX\,CY}{CA\,CB}=-1$$
hods  if and only if $X$ and $Y$ are isogonally conjugate of type II
relative to $A$}.

\vskip 2mm
In Section 2 we give an exact description of the geometric configurations,
which are characterized by the above equalities.

The aim of our considerations is to apply the above mentioned problems in order to
clear up the geometric meaning of the solution to the classical Fermat problem
for a triangle and arbitrary weights (cf. \cite{JK}).

In Section 3 we consider the Fermat problem with positive weights:
\vskip 2mm {\it Given a $\triangle ABC$ and a point $Y$ in its
plane. If $\lambda, \mu, \nu$ are three positive numbers, describe
geometrically the point $Y$, which minimizes the function}
$$\lambda\,AY+\mu\,BY+\nu\,CY.$$

We prove that in the general case the given numbers determine uniquely a point
$X$ and the solution to the problem is the point $Y$, isogonally conjugate of type $I$
to the point $X$.
\vskip 2mm
In Section 4 we consider the Fermat problem for two positive and one negative weights:
\vskip 2mm
{\it Given a $\triangle ABC$ and a point $Y$ in its plane. If $\lambda, \mu, \nu$ are
three positive numbers, describe geometrically the point $Y$, which minimizes the
function}
$$-\lambda\,AY+\mu\,BY+\nu\,CY.$$

We prove that in the general case the given numbers determine uniquely a point $X$ and
the solution to the problem is the point $Y$, isogonally conjugate of type II relative
to $A$ to the point $X$.

\section{Inequalities and equalities, characterizing the isogonal conjugacy in Euclidean plane}

Let be given a $\triangle ABC$ and its circumscribed circle $k$.
Denote by $\imath$ the map of the isogonal conjugacy with respect to the
given triangle. We shall exclude the points on the circle $k$,
different from the vertices of the triangle because we do not
consider points at infinity.

There are shown below (Figure 1) the regions in the angular domains
$\angle DAE$ and $\angle FAG$, determined by the map $\imath$, the
given $\triangle ABC$, the circle $k$ and the point $A$.
\vskip 2mm
\begin{center}
\includegraphics[width=5cm]{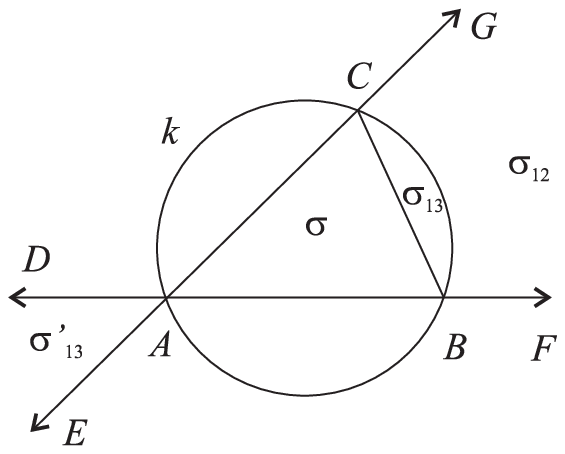}
\vskip 2mm
Figure 1
\end{center}
\vskip 2mm
According to the definition of $\imath$ we have:

1) $\imath(\sigma)=\sigma$.

If $M$ is a boundary point for the domain $\sigma$, then: $\imath
(M)=A$ for any point $M$ on the segment $BC$; \; $\imath (M)=B$ for
any point $M$ on the segment $CA$; \; $\imath (M)=C$ for any point
$M$ on the segment $AB$.

2) $\imath(\sigma_{12})=\sigma_{12}$.

If $M$ is a boundary point for the domain $\sigma_{12}$, then:
$\imath (M)=C$ for any point $M$ on the ray $BF^{\rightarrow}$; \;
$\imath (M)=B$ for any point $M$ on the ray $CG^{\rightarrow}$.

3) $\imath(\sigma_{13})=\sigma'_{13}$, \; $\imath(\sigma'_{13})=\sigma_{13}$.

If $M$ is a boundary point for the domain $\sigma_{13}\cup
\sigma'_{13}$, then: $\imath (M)=A$ for any point $M$ on the segment
$BC$; \; $\imath(M)=C$ for any point $M$ on the ray
$AD^{\rightarrow}$; \; $\imath(M)=B$ for any point $M$ on the ray
$AE^{\rightarrow}$.

Further we shall make difference between the cases, when the points are in the domains
$\sigma, \, \sigma_{12}$, or $\sigma_{13}\cup \sigma'_{13}$. We give the following definitions:

Two isogonally conjugate points $X$ and $Y$ are said to be:

1) \emph{isogonally conjugate of type $I$}, if $X,\,Y$ either are in $\sigma$ or
are boundary points for $\sigma$;

2) \emph{isogonally conjugate of type $II$ relative to $A$}, if $X,\,Y$
either are in $\sigma_{12}$ or are boundary points for
$\sigma_{12}$.

3) \emph{isogonally conjugate of type $III$ relative to $A$}, if $X,\,Y$
either are in $\sigma_{13}\cup \sigma'_{13}$ or are boundary points
for $\sigma_{13}\cup \sigma'_{13}$.
\vskip 2mm

Given a $\triangle ABC$ and two arbitrary points $X$, $Y$ in its plane.
First we shall consider the following problem.
\vskip 2mm
{\bf Problem 1.} {\it Prove that
$$\frac{AX AY}{AB AC}+\frac{BX BY}{BA BC}+\frac{CX CY}{CA CB}\geq 1, \leqno(1)$$
and the equality holds  if and only if $X$ and $Y$ are isogonally
conjugate of type $I$.} \vskip 2mm {\bf Solution:} In this section,
for convenience in the computations we identify the Euclidean plane
of $\triangle ABC$ with the Gauss plane of complex numbers. As
usual, the complex numbers (affixes) corresponding to the points,
e.g. $A,B,C$, are denoted by small letters: $a,b,c$ etc. Then the
required inequality can be written in the form
$$\left | \frac{(x-a)(y-a)}{(b-a)(c-a)}\right |
+\left | \frac{(x-b)(y-b)}{(a-b)(c-b)}\right | +\left
|\frac{(x-c)(y-c)}{(a-c)(b-c)}\right |\geq 1.\leqno(2)$$

The last inequality follows from the standard inequality
$$\left | \frac{(x-a)(y-a)}{(b-a)(c-a)}\right |
+\left | \frac{(x-b)(y-b)}{(a-b)(c-b)}\right | +\left
|\frac{(x-c)(y-c)}{(a-c)(b-c)}\right |\geq$$
$$
\left | \frac{(x-a)(y-a)}{(b-a)(c-a)} +\frac{(x-b)(y-b)}{(a-b)(c-b)}
+\frac{(x-c)(y-c)}{(a-c)(b-c)}\right |$$ and the fact that
$$\frac{(x-a)(y-a)}{(b-a)(c-a)}
+\frac{(x-b)(y-b)}{(a-b)(c-b)}
+\frac{(x-c)(y-c)}{(a-c)(b-c)}=1.\leqno(3)$$

The equality in (2) holds  if and only if the three complex numbers
$$\frac{(x-a)(y-a)}{(b-a)(c-a)}, \quad \frac{(x-b)(y-b)}{(a-b)(c-b)},
\quad \frac{(x-c)(y-c)}{(a-c)(b-c)}\leqno{(4)}$$ are real and
non-negative. The complex numbers (4) are real if and only the
points $X$ and $Y$ are isogonally conjugate with respect to the $\triangle ABC$.
Since at most one of the numbers (4) can be zero, then we have
the following cases:

i) The three numbers (4) are positive. This case occurs if and only if
$X$ and $Y$ are interior points for $\triangle ABC$ and isogonally conjugate.

ii) One of the numbers (4) is zero. This case occurs if and only if one of the points
$X$ and $Y$ is a vertex of the triangle and the other is an interior point for the
opposite side. This is the case when the points $X$ and $Y$ are boundary for
$\triangle ABC$ and isogonally conjugate.

Thus we obtained:

{\it The equality in $(1)$ holds  if and only if the points $X$ and
$Y$ are isogonally conjugate of type $I$.} \hfill{\qed}

\vskip 2mm

If $BC=a, \; CA=b$ è $AB=c$, then the condition that two points $X$ and
$Y$ are isogonally conjugate of type $I$ is characterized
analytically by the equality
$$a\,AX\,AY + b\,BX\,BY + c\,CX\,CY=abc.$$
\vskip 2mm
{\bf Problem 2.} {\it Prove that
$$-\frac{AX AY}{AB AC}+\frac{BX BY}{BA BC}+\frac{CX CY}{CA CB}\geq -1, \leqno(5)$$
and the equality holds if and only if the points $X$ and $Y$ are
isogonally conjugate of type $II$ relative to $A$.}
\vskip 2mm
{\bf Solution:}
Similarly to the solution of Problem 1 the inequality (5) can be written in the form
$$\left | \frac{(x-a)(y-a)}{(b-a)(c-a)}\right |
-\left | \frac{(x-b)(y-b)}{(a-b)(c-b)}\right | -\left
|\frac{(x-c)(y-c)}{(a-c)(b-c)}\right |\leq 1.\leqno(6)$$

The last inequality follows from the standard inequality
$$\left | \frac{(x-a)(y-a)}{(b-a)(c-a)}\right |
-\left | \frac{(x-b)(y-b)}{(a-b)(c-b)}\right | -\left
|\frac{(x-c)(y-c)}{(a-c)(b-c)}\right |\leq$$
$$
\left | \frac{(x-a)(y-a)}{(b-a)(c-a)} +\frac{(x-b)(y-b)}{(a-b)(c-b)}
+\frac{(x-c)(y-c)}{(a-c)(b-c)}\right |$$ taking into account (3).

The equality in (6) holds  if and only if
$$\frac{(x-a)(y-a)}{(b-a)(c-a)}\geq 0, \quad \frac{(x-b)(y-b)}{(a-b)(c-b)}\leq 0,
\quad \frac{(x-c)(y-c)}{(a-c)(b-c)}\leq 0. \leqno (7)$$

The above complex numbers are real if and only if the points $X$ and
$Y$ are isogonally conjugate with respect to $\triangle ABC$. Further, as in
the solution of the Problem 1 we have:

i) The three numbers (7) are different from zero if and only if the points $X$ and $Y$
are in the domain $\sigma_{12}$ and are isogonally conjugate.

ii) One of the complex numbers (7) is zero. This is the case when the points $X$ and $Y$
are boundary for the domain $\sigma_{12}$ and isogonally conjugate.

Thus we proved:

{\it The equality in $(5)$ holds  if and only if the points $X$ and
$Y$ are isogonally conjugate of type $II$ relative to $A$.} \hfill{\qed}
\vskip 2mm

The condition two points $X$ and $Y$ to be isogonally conjugate of type $II$ relative
to $A$ is characterized analytically by the equality
$$-a\,AX\,AY + b\,BX\,BY + c\,CX\,CY=-abc.$$

Thus the he following natural question arises: Does there exist an
inequality, similar to the inequalities in the problems 1 and 2,
characterizing the isogonal conjugacy of type $III$?

In this case the answer is negative, but one can prove the following
property of the isogonal conjugacy with respect to $\triangle ABC$:
\vskip 2mm
{\bf Problem 3.} {\it Given two pints $X$ and $Y$ isogonally
conjugate with respect to $\triangle ABC$. Prove that $X$ and $Y$ are
isogonally conjugate of type $III$ relative to $A$ if and only if}
$$-\frac{AX AY}{AB AC}+\frac{BX BY}{BA BC}+\frac{CX CY}{CA CB}=1, \leqno(8)$$
\vskip 2mm
{\bf Solution:} Since $X$ and $Y$ are isogonally conjugate, then the three complex numbers
in (4) are real.

If $X$ and $Y$ are isogonally conjugate of type $III$ relative to $A$,
then the first number in (4) is non-positive and the other two are
non-negative. Then
$$-\left | \frac{(x-a)(y-a)}{(b-a)(c-a)}\right |
+\left | \frac{(x-b)(y-b)}{(a-b)(c-b)}\right | +\left
|\frac{(x-c)(y-c)}{(a-c)(b-c)}\right |$$
$$=\frac{(x-a)(y-a)}{(b-a)(c-a)}+\frac{(x-b)(y-b)}{(a-b)(c-b)}+
\frac{(x-c)(y-c)}{(a-c)(b-c)}=1.$$

Conversely, let
$$-\left | \frac{(x-a)(y-a)}{(b-a)(c-a)}\right |
+\left | \frac{(x-b)(y-b)}{(a-b)(c-b)}\right | +\left
|\frac{(x-c)(y-c)}{(a-c)(b-c)}\right |=1,$$

Taking into account (3) and the arguments in the solutions of the
problems 1 and 2, it follows that the first number in (4) is
non-positive and the other two are non-negative. This means that the
points $X$ and $Y$ are isogonally conjugate of type $III$ relative to
$A$. \hfill{\qed}
\vskip 2mm
Thus we obtained:

{\it Two isogonally conjugate points $X$ and $Y$ are isogonally
conjugate of type $III$ relative to $A$ if and only if}
$$-a\,AX\,AY + b\,BX\,BY + c\,CX\,CY=abc.$$

\section{Applications to the Fermat problem with positive weights}

For a $\triangle ABC$ we use the standard denotations: $BC=a, \, CA=b,
\, AB=c; \; \angle A=\alpha, \, \angle B=\beta, \, \angle C=\gamma;
\; k(O,R)$ - the circumscribed circle of $\triangle ABC$; $S$ - the
area of $\triangle ABC$.

Let $M\not\in k$ be an arbitrary point and $\varphi (M, r)$ is an
inversion with center $M$ and an arbitrary radius $r$. If $\varphi
(A)=A',\, \varphi(B)=B',\, \varphi(C)=C'$, then
$$B'C':C'A':A'B'=a \, AM: b \, BM: c \, CM.$$

Let us denote by $A_1B_1C_1$ the pedal triangle of the point $M$
with respect to the given $\triangle ABC$ and let its angles be $\angle
A_1=\alpha_1,\, \angle B_1=\beta_1,\,\angle C_1=\gamma_1$,
respectively. We have the following formulas for the sides of the
pedal triangle
$$B_1C_1=\frac{a}{2R}AM,\quad C_1A_1=\frac{b}{2R}BM,\quad A_1B_1=\frac{c}{2R}\,CM.$$
Hence $\triangle A_1B_1C_1 \sim \triangle A'B'C'$.

Thus, the inverse question arises: How to find the point $M$, if the
angles of $\triangle A'B'C'$ \; $(\triangle A_1B_1C_1)$ are known?

We shall use the following statement.
\vskip 2mm
{\bf Theorem.} \cite{G} {\it Let $(\alpha_1, \, \beta_1, \,\gamma_1)$ be three angles,
which are angles of a triangle. Then:

i) There exists a unique point $M$, which is an interior point for the circle $k$, such that
the angles of its pedal $\triangle A_1B_1C_1$ are $(\alpha_1, \, \beta_1, \,\gamma_1)$,
respectively. In this case $\triangle A_1B_1C_1$ has the orientation of the basic $\triangle ABC$.

ii) If  $(\alpha_1, \, \beta_1, \,\gamma_1)\neq
(\alpha,\,\beta,\,\gamma)$, then there exists a unique point $N$,
which is an exterior point for the circle $k$, such that the angles
of its pedal $\triangle A_2B_2C_2$ are $(\alpha_1, \, \beta_1,
\,\gamma_1)$, respectively. In this case the orientation of
$\triangle A_2B_2C_2$ is opposite to the orientation of the basic
$\triangle ABC$.

iii) The points $M$ and $N$ are inversive with respect to the circle $k$.}
\vskip 3mm
Further we give the relations between the angles $(\alpha,\,\beta,\,\gamma)$ and
$(\alpha_1, \, \beta_1, \,\gamma_1)$ which characterize the condition $M$ to belong to
$\sigma, \, \sigma_{13}, \, \sigma_{12}$ or $\sigma'_{13}$.

Let $M$ be an \emph{interior} point for the circumscribed circle $k$. Then we have (Figure 2, a):
$$\angle BMC=\alpha+\alpha_1,\quad \angle
CMA=\beta+\beta_1,\quad \angle AMB=\gamma+\gamma_1. \leqno(9)$$

Here we adopt the convention that $\angle BMC>\pi$ if and only if the points $A$ and $M$
lie on different sides of the chord $BC$ (Figure 2, b) and $\angle BMC=\pi$ if and only if
$M$ lies on the chord $BC$. Then the formulas (9) are again valid.

$$\begin{array}{cc}
\includegraphics[width=5cm]{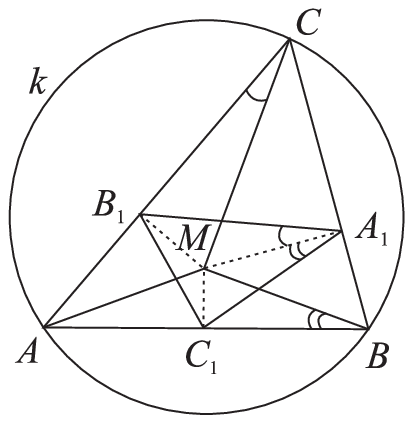} \qquad & \qquad \includegraphics[width=5cm]{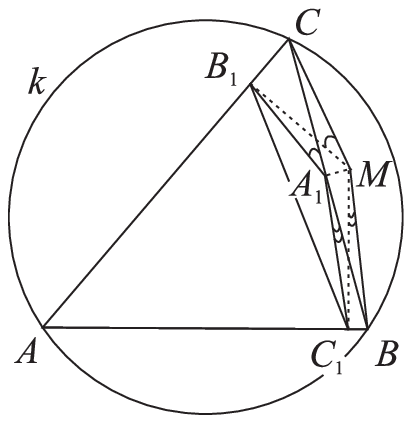}\\
[2mm]
{\rm Figure\; 2,\; a} \qquad & \qquad {\rm Figure \;2,\; b}
\end{array}$$

Thus we have:
$$M \in \sigma \iff \alpha+\alpha_1<\pi, \quad \beta + \beta_1<\pi, \quad \gamma+\gamma_1<\pi;$$
$$M \in (BC) \iff \alpha+\alpha_1=\pi, \quad \beta + \beta_1<\pi, \quad \gamma+\gamma_1<\pi;$$
$$M \in \sigma_{13} \iff \alpha+\alpha_1>\pi, \quad \beta + \beta_1<\pi, \quad \gamma+\gamma_1<\pi;$$

Formulas (9) give a natural way to construct the point $M$.

If $M$ is an interior point for the domain $\sigma$ and $N$ is the isogonally conjugate
point to $M$, then the equality $\angle BMC+\angle BNC=\pi+\alpha$ and (9) imply that
$$\angle BNC=\pi-\alpha_1,\quad \angle CNA=\pi-\beta_1,\quad
\angle ANB=\pi-\gamma_1. \leqno(10)$$
\vskip 2mm
Let $M$ be an \emph{exterior} point for the circumscribed circle $k$.

If the points $M$ and $A$ lie on opposite sides of the line $BC$ (Figure 3), then
$\angle BMC=\alpha_1-\alpha>0$.
\vskip 2mm
\begin{center}
\includegraphics[width=5cm]{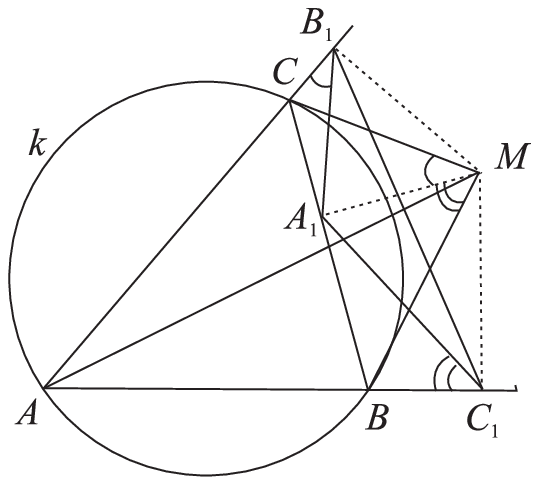}
\vskip 2mm
Figure 3
\end{center}
\vskip 2mm
If the points $M$ and $A$ lie on the same side of the line $BC$, then it follows
in a similar way that $\angle BMC=\alpha-\alpha_1>0$.

Therefore the domains $\sigma_{12}$ and $\sigma'_{13}$ are characterized as follows:
$$M\in \sigma_{12} \iff \alpha < \alpha_1, \quad \beta>\beta_1, \quad \gamma>\gamma_1;$$
$$M\in \sigma'_{13} \iff \alpha > \alpha_1, \quad \beta < \beta_1, \quad \gamma < \gamma_1.$$

Now we shall consider the Fermat problem with positive weights. Using the inequalities
from Problem 1 we shall give a geometric interpretation of the solution to the Fermat
problem.
\vskip 2mm
{\bf Problem 4.} {\it Given a $\triangle ABC$ and three positive numbers
$\lambda, \mu, \nu$. Find a point $Y$ in the plane of the given triangle, which
minimizes the function}
$$F(Y)=\lambda \,AY+\mu \,BY + \nu \,CY \leqno(11)$$
\vskip 2mm
{\bf Solution:} 1) First, let one of the given numbers $\lambda, \mu, \nu$ is greater
or equal to the sum of the other two, for example, $\lambda \geq \mu + \nu$. Then for
any point $Y$ we have
$$\begin{array}{l}
F(Y) \geq (\mu + \nu) AY+\mu BY + \nu\, CY\\
[2mm]
=\mu(AY+BY)+\nu(AY + CY) \geq  F(A).\end{array}$$

The equality holds  if and only if $Y\equiv A$. \vskip 2mm 2) Let
now the given numbers satisfy the triangle inequality:
$$\lambda+\mu>\nu, \quad \mu +\nu>\lambda,\quad\nu+\lambda>\mu.$$
Then there exists a triangle, whose sides are $\lambda, \mu, \nu$.
Denote the angles of this triangle by $(\alpha_1, \beta_1, \gamma_1)$, respectively,
and consider the function
$$f(Y)=\sin \alpha_1 \, AY + \sin \beta_1 \, BY + \sin \gamma_1 \, CY.$$

The investigation of the function $f(Y)$ is equivalent to the study of (11).

According to Theorem 1, there exists a uniquely determined point
$X$, which is an interior point for the circumscribed circle $k$,
such that the angles of its pedal $\triangle A_1B_1C_1$ are $\angle
A_1=\alpha_1, \, \angle B_1=\beta_1, \, \angle C_1=\gamma_1$,
respectively. If $R_1$ is the circumradius of $\triangle A_1B_1C_1$,
then the expression for $f(Y)$ can be written in the form
$$f(Y)=\frac{1}{4RR_1}(a\,AX\,AY+b\,BX\,BY+c\,CX\,CY).\leqno(12)$$

We shall consider the following cases:

$\bullet$ the point $X$ is interior or boundary for $\triangle ABC$;

$\bullet$ the point $X$ is exterior for $\triangle ABC$. \vskip 2mm

In the first case, applying the statement of Problem 1 to (12), we
obtain that the minimum of $f(Y)$ is $\displaystyle{\frac{S}{R_1}}$
and the minimum is attained if and only if $Y$ is isogonally
conjugate to the point $X$.

More precisely, $X$ is an interior point for $\triangle ABC$ if and
only if $\alpha+\alpha_1<\pi, \beta+\beta_1<\pi,
\gamma+\gamma_1<\pi$. In this case $f_{\rm
min}=\displaystyle{\frac{S}{R_1}}$ is attained for the point $Y$,
which is isogonally conjugate to $X$  and it is also an interior
point for $\triangle ABC$.

The point $X$ lies on the side $BC$ of $\triangle ABC$ if and only
if $\alpha+\alpha_1=\pi$. In this case $f_{\rm
min}=\displaystyle{\frac{S}{R_1}}$ is attained for the vertex $Y=A$,
which is isogonally conjugate to $X$. \vskip 2mm

Let now $X$ be an exterior point for $\triangle ABC$ (Figure 4). Then we can
assume that $\alpha + \alpha_1>\pi$.
\vskip 2mm
\begin{center}
\includegraphics[width=5cm]{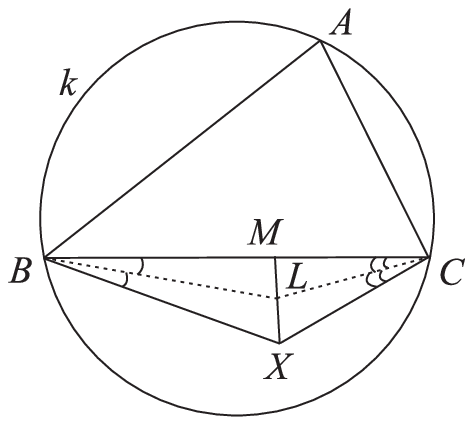}
\vskip 2mm
Figure 4
\end{center}
\vskip 2mm
\noindent
Let $XM$ be the bisector of $\angle BXC$ in $\triangle
BCX$. Then there exists a number $q>1$ such that $$BX=q BM, \quad CX=q CM.$$

Let us consider the Apollonian circle $k_0$ with basic points $M, \,
X$ and a ratio $\displaystyle{\frac{BX}{BM}=\frac{CX}{CM}=q}$. If
$L$ is the intersection point of the bisectors of the angles
$\angle XBM$ and $\angle XCM$ with the segment $MX$, then $k_0$
passes through the points $B,C,L$ and consequently $A$ is an
interior point for $k_0$. Hence $\displaystyle{\frac{AX}{AM}>
\frac{BX}{BM}=q.}$

Replacing $BX=q BM$ and $CX=q CM$ in (12), we find
$$f(Y)=\frac{1}{4RR_1}\{q(a \,AM\,AY+b\,BM\,BY+c\,CM\,CY)+a(AX-q\,AM)AY\}.$$
Since $AX-q\,AM>0$, then it follows from Problem 1 that the minimum
of $f(Y)$ is $\displaystyle{\frac{qS}{R_1}}$ and it is attained if
and only if $Y$ is isogonally conjugate to the point $M$, i.e.
$Y\equiv A$.

$\hfill{\qed}$

A natural way for constructing the minimizing point $Y$ in the general case
follows from equalities (10).

\section{Applications to the Fermat problem with one negative and two positive weights}

In this section we consider the Fermat problem for a triangle with
mixed (one negative and two positive) weights. We show that the
solution to this problem can be obtained from the solution of the
Fermat problem with positive weights.

{\bf Problem 5.} {\it Given a $\triangle ABC$ and three positive
numbers $\lambda, \mu, \nu$. Find a point $Y$ in the plane of the
triangle, which minimizes the function}
$$G(Y)=-\lambda \,AY+\mu \,BY + \nu \,CY.\leqno(13)$$
\vskip 2mm {\bf Solution:} We shall use the following relation
between the problems 4 and 5. \vskip 2mm {\bf Lemma.} \cite{JK} {\it
If $Q$ is a solution of Problem $5$, then:

i) $Q$ does not lie in the same half-plane with respect to the line $BC$
as the point $A$;

ii) For any point $D$ lying on the ray, which is the opposite to the
ray $QA^{\rightarrow}$, the point $Q$ is a solution of Problem $4$
for $\triangle BCD$ and the positive numbers $\lambda, \mu, \nu$.}
\vskip 2mm
{\bf Proof:} i) The statement follows from the obvious
fact that if the point $Q$ lies in the half-plane with respect to the line
$BC$, containing $A,$ and $Q'$ is the symmetric point to $Q$ with respect to
this line, then $G(Q')<G(Q)$ (Figure 5, a).

ii) Let $Y$ be an arbitrary point in the plane of $\triangle ABC$ and
$$F_D(Y)=\lambda \,DY+\mu \,BY + \nu \,CY.$$
Then
$$F_D(Y)-F_D(Q)-(G(Y)-G(Q))=\lambda(DY+YA-(DQ+QA))=\lambda(DY+YA-DA)\ge 0.$$
If $Q$ is a solution of Problem $5$, it follows that
$$F_D(Y)-F_D(Q)\ge(G(Y)-G(Q))\ge 0,$$
and consequently $Q$ is a solution of Problem $4$ for $\triangle
BCD$ and the positive weights $\lambda, \mu, \nu$ (Figure 5, b).
\hfill{\qed} \vskip 2mm
$$\begin{array}{cc}
\includegraphics[width=4cm]{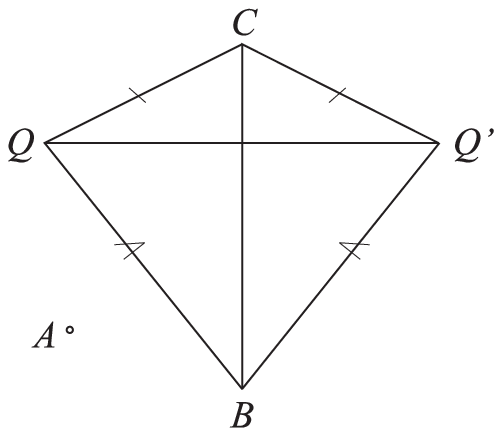} \qquad & \qquad \includegraphics[width=4cm]{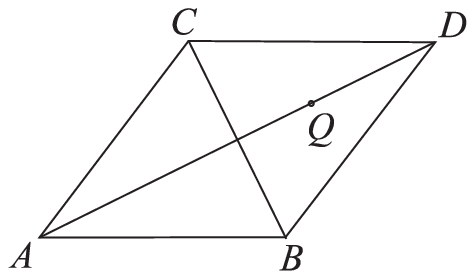}\\
[2mm]
 {\rm Figure \; 5 \; a} \qquad & \qquad {\rm Figure \; 5 \; b}
\end{array}$$
\vskip 2mm The lemma implies that any solution $Q$ of Problem 5 (if
it exists) is: $B$, $C$ or a point in the domain $\sigma_{13} \cup \,
{\rm arc} \, AB \cup \sigma_{12}$.

First we note that if $\lambda>\mu+\nu,$ then $$G(Y)\leq (\mu+\nu-\lambda)AY +\mu AB+\nu AC$$
and therefore \emph{there is no solution} to the problem.

Now let $\lambda\le\mu+\nu.$ Then $G(Y)\geq -G(A)$ and consequently
Problem 5 {\it has} a solution. It follows from the lemma and
Problem 4 that the solution is: \vskip 1mm {\it The point $B,$ if
$\mu>\lambda+\nu; $ the \, point $C,$ if $\nu>\lambda+\nu;$ the \,
point $B$ and / or $C,$ if $\lambda=\mu+\nu.$} \vskip 1mm

Now we shall consider the case when the numbers $\lambda, \mu,\nu$
satisfy the triangle inequalities:
$$\lambda+\mu>\nu, \quad \mu+\nu>\lambda, \quad \nu+\lambda>\mu.$$

Let us denote the angles of the triangle, whose sides are $\lambda, \mu, \nu$ with
$\alpha_1, \beta_1, \gamma_1$, respectively, and consider the function
$$g(Y)=-\sin \alpha_1 AY+\sin \beta_1 BY+\sin \gamma_1 CY.\leqno (14)$$

The investigation of the function $g(Y)$ is equivalent to the study of (13).

Assume that the point $Q$ from the lemma lies in the domain $\sigma_{13}$
(Figure 6).
$$\begin{array}{cc}
{\rm Figure \; 6} \qquad & \qquad \includegraphics[width=5cm]{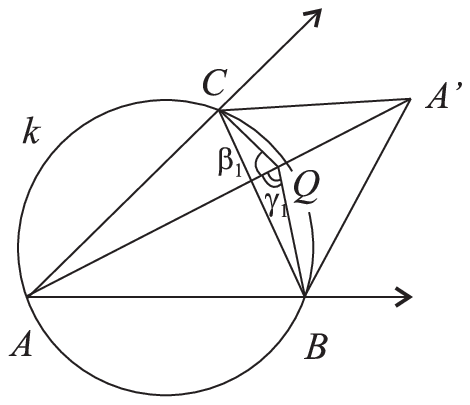}
\end{array}$$
\vskip 2mm
Let $\triangle BCA'$ satisfy the conditions of the lemma. Then it follows from
the solution of Problem 4 that $\angle BQA'=\pi-\gamma_1, \angle CQA'=\pi-\beta_1$ and
$\angle AQB=\gamma_1, \angle AQC=\beta_1$.

Moreover, we have
$$\beta_1>\beta, \quad \gamma_1>\gamma. \leqno(15)$$

Let us also consider the point $P$, exterior for the circle $k$,
whose pedal triangle has angles $\alpha_1, \beta_1, \gamma_1$,
respectively. Because of (15), the point $P$ lies in the domain
$\sigma'_{13}$, and because of the equalities $\angle
BQC=\pi-\alpha_1, \angle AQC=\beta_1, \angle AQB=\gamma_1$, the
points $P$ and $Q$ are isogonally conjugate of type III relative to $A$.
Then we have from Problem 3 that
$$g(Q)=\frac{abc}{4RR_1},$$
where $R_1$ is the circumradius of the pedal triangle of the point
$P$.

On the other hand, we find the following expression for the value of the function $g$ in
the point $B$:
$$g(B)= \frac{ac}{4RR_1}(PC-PA)<\frac{abc}{4RR_1}=g(Q),$$
which contradicts to the property of $Q$.

Hence $Q$ can not be in the domain $\sigma_{13}$.

Let $Q$ lie on the arc $BC$ of $k$, which does not contain the point $A$.
Then $\beta_1=\beta, \gamma_1=\gamma$ and $\alpha_1=\alpha$.
Therefore
$$g(Y)=\frac{1}{2R}(-a AY+b BY+c CY)\geq 0,$$
and the equality is attained for the points of the arc $BC$.

In this case, a solution of the problem is any point on the arc
$BC$, which does not contain $A$, and $g_{\rm min}=0$ (Ptolemy's theorem).

Let $Q \in \sigma_{12}$. Then $\beta_1<\beta, \gamma_1<\gamma$.

Let us consider the point $P$, exterior for $k$, whose pedal
triangle has angles $\alpha_1, \beta_1, \gamma_1$. As in the above
we conclude that the points $P$ and $Q$ are isogonally conjugate of
type $II$ relative to $A$.

In this case the solution is the point $Q$, which is isogonally conjugate to the point
$P$ and
$$g_{\rm min}=-\frac{abc}{4RR_1}=-\frac{S}{R_1}.$$

In all other cases the solution is the point $B$, the point $C$ or
both of them in accordance with $g(B)<g(C)$, $g(C)<g(B)$ or
$g(B)=g(C),$ respectively.

Thus we obtained:

{\it $1.$ If $\beta_1=\beta, \gamma_1=\gamma$, then the solutions of
the problem are exactly the points on the arc $BC$, not containing the point $A$,
and $g_{\rm min}=0$.

$2.$ If $\beta_1<\beta, \gamma_1<\gamma$, then the solution is the point $Q$, which
is isogonally conjugate of type II relative to $A$ to the point $P$, and
$$g_{\rm min}=-\frac{abc}{4RR_1}=-\frac{S}{R_1}.$$

$3.$ Let $\beta_1>\beta$ or $\gamma_1>\gamma$.
\vskip 2mm
\hskip 10mm
$3.1.$  If \, $\displaystyle{\frac{\sin \beta_1-\sin \gamma_1}{\sin \alpha_1}
>\frac{\sin \beta-\sin \gamma}{\sin \alpha}}$,
then the solution is the point $B$;
\vskip 2mm
\hskip 10mm
$3.2.$
If \, $\ds{\frac{\sin \beta_1-\sin \gamma_1}{\sin \alpha_1}
=\frac{\sin \beta-\sin \gamma}{\sin \alpha}}$, then the
solutions are the point $B$ and the point $C$;
\vskip 2mm
\hskip 10mm
$3.2$  If \, $\ds{\frac{\sin \beta_1-\sin \gamma_1}{\sin \alpha_1}
<\frac{\sin \beta-\sin \gamma}{\sin \alpha}}$, then the solution is the point $C$}.
\vskip 2mm
{\bf Remark.} Standard computations show that
$${\frac{\sin \beta_1-\sin \gamma_1}{\sin \alpha_1}
>\frac{\sin \beta-\sin \gamma}{\sin\alpha}}\iff
\tan\frac{\beta_1}{2}\cot\frac{\gamma_1}{2}>\tan\frac{\beta}{2}\cot\frac{\gamma}{2}.$$
Hence, the solution in the domain $\alpha_1\neq\alpha,$ $\beta_1\ge\beta,$
$\gamma_1\le\gamma$ ($\alpha_1\neq\alpha,$ $\beta_1\le\beta,$
$\gamma_1\ge\gamma$) is the point $B$ ($C$).

\end{document}